\documentclass[11pt]{article}
\usepackage{multicol}
\usepackage{amssymb,amsmath,amsthm,bm}
\usepackage[colorlinks=true, urlcolor=blue,linkcolor=blue, citecolor=blue]{hyperref}
\usepackage{mathrsfs,verbatim}
\usepackage{caption,cleveref }
\usepackage{graphicx, float}
\usepackage{color, mathtools, tikz, xcolor}
\usepackage{enumerate,enumitem}
\usepackage{tikz}
\usetikzlibrary{shapes.geometric}
\usepackage{bbm}

\usepackage[UKenglish]{babel}
\usepackage[UKenglish]{isodate}

\usetikzlibrary{calc}

\usepackage[mathscr]{euscript}
\usepackage{appendix}
\usepackage{geometry}
\numberwithin{equation}{section}

	\newtheorem{theorem}{Theorem}[section]

	\newtheorem{corollary}[theorem]{Corollary}

	\newtheorem{proposition}[theorem]{Proposition}
	\newtheorem{lemma}[theorem]{Lemma}

\def\al#1{}
	\renewcommand{\al}[1]{\footnote{\textbf{AL: }#1}} 

\def\eps{{\varepsilon}}

\title{From finding a spanning subgraph $H$ to an $H$-factor}

\author{Allan Lo\thanks{School of Mathematics, University of Birmingham, B15 2TT, United Kingdom. Email: s.a.lo@bham.ac.uk.
Allan Lo was partially supported by EPSRC, grant no. EP/V048287/1.
}
}
\date{\today}

\begin{document}

\maketitle

\begin{abstract}
A typical Dirac-type problem in extremal graph theory is to determine the minimum degree threshold for a graph~$G$ to have a spanning subgraph $H$, e.g. the Dirac theorem. 
A natural follow-up problem is to seek an $H$-factor, which is a spanning set of vertex-disjoint copies of~$H$.
In this short note, we present a method for obtaining an upper bound on the minimum degree threshold for an $H$-factor from one for finding a spanning copy of~$H$. 

As an application, we prove that, for all $\eps>0$ and sufficiently large~$\ell$, any oriented graph~$G$ on $\ell m$ vertices with minimum semi-degree $\delta^0(G) \ge  (3/8+ \eps)\ell m$ contains a $C_\ell$-factor, where $C_\ell$ is an arbitrary orientation of a cycle on $\ell$ vertices. 
This improves a result of Wang, Yan and Zhang. 
\end{abstract}

\section{Introduction}

A fundamental problem in extremal graph theory is to determine the minimum degree threshold for the existence of a spanning subgraph. 
The celebrated theorem of Dirac~\cite{Dirac} states that any graph on $n \ge 3$ vertices with minimum degree at least $n/2$ contains a Hamilton cycle (which is a spanning cycle). 
This theorem has been generalised to various settings such as digraphs and hypergraphs, see surveys~\cite{MR2889513,MR3526407}, respectively. 

Given a graph $H$, an \emph{$H$-tiling} is a set of vertex-disjoint copies of~$H$. 
We say that a graph~$G$ contains an \emph{$H$-factor} (or a \emph{perfect $H$-tiling}) if $G$ contains a spanning $H$-tiling. 
Clearly, $|H|$ divides~$|G|$ is a necessary condition for an $H$-factor in~$G$, which we will assume holds whenever we talk about $H$-factors.
In addition, we will assume that $n=|G|$ is sufficiently large. 
A classical result of Hajnal and Szemer\'edi~\cite{HajnalSzemeredi} states that $\delta(G) \ge (t-1)n/t$ implies a $K_t$-factor. 
For a graph~$H$, K\"uhn and Osthus~\cite{MR2506388} determined the minimum degree threshold for an $H$-factor up to an additive constant. 
Tiling problems have also been considered for digraphs and hypergraphs, see surveys~\cite{MR2588541,MR3526407}, respectively. 

The main contribution of this paper is to provide a simple method of converting the minimum degree threshold for the existence of a spanning graph~$H$ into a (decent) minimum degree threshold for the existence of an $H$-factor.
We will illustrate it via finding cycle-factors in oriented graphs.

\subsection{Hamilton cycles in oriented graphs}

A Hamilton cycle in a digraph is \emph{directed} if all its edges are oriented in the same direction. 
Recall that the minimum semi-degree of a digraph~$G$ denoted by $\delta^0(G)$ is the minimum of the in-degree and out-degree over all vertices, that is, $\delta^0(G) = \min \{d^-(v), d^+(v): v \in V(G)\}$.
The oriented graph analogue of Dirac's theorem was proved asymptotically by Kelly, K\"uhn and Osthus~\cite{MR2454564} and then exactly by Keevash, K\"uhn and Osthus~\cite{MR2472138}. 
Instead of a directed Hamilton cycle, one may ask for an arbitrary orientation of a Hamilton cycle.
Kelly~\cite{Kelly_arborientation} gave an asymptotically sharp threshold and, recently, Wang, Wang and Zhang~\cite{wang2025arbitrary} determined the exact threshold. 
\begin{theorem}[Wang, Wang and Zhang~\cite{wang2025arbitrary}] \label{thm:arboriHC}
There exists an integer $n_0$ such that any oriented graph~$G$ on $n\ge n_0$ vertices with minimum semi-degree $\delta^0(G) \ge \lceil (3n-1)/8 \rceil$ contains every possible orientation of a Hamilton cycle.
\end{theorem}

\subsection{Cycle-factors in oriented graphs}

For $\ell \in \mathbb{N}$, let $\overset{\rightarrow}{C}_\ell$ be the directed cycle on $\ell$ vertices. 
When $\ell =3$, Keevash and Sudakov~\cite{KeevashSudakov} found an oriented graph~$G$ with $\delta^0(G) = (n-1)/2-1$ and $n \equiv 3 \mod{18}$ without a $\overset{\rightarrow}{C}_{3}$-factor.
Molla and Li~\cite{MR4025428} proved that if $G$ does not have a `divisibility barrier', then $\delta^0(G) \ge (1/2-c)n$ for some $c >0$ is already sufficient for a $\overset{\rightarrow}{C}_{3}$-factor.
In particular, they showed that every regular tournament~$G$ on $n$ vertices, that is when $\delta^0(G) = (n-1)/2$, contains a $\overset{\rightarrow}{C}_{3}$-factor.
When $\ell \ge 4$, Wang, Yan and Zhang~\cite{MR4763245} showed that $\delta^0(G) \ge (1/2-c) n$ implies a $\overset{\rightarrow}{C}_{\ell}$-factor. 
We improve their result for large~$\ell$.

\begin{theorem} \label{thm:arboriHC-tiling}
For all $\eps >0$, there exists $\ell_0 = \ell_0(\eps)$ such that for all $\ell \ge \ell_0$ and any oriented cycle~$C_\ell$ on $\ell$ vertices, any oriented graph~$G$ on $ \ell m$ vertices with minimum semi-degree $\delta^0(G) \ge  (3/8+ \eps)\ell m$ contains $C_\ell$-factor.
\end{theorem}

On the other hand, the constant $3/8$ cannot be reduced. 
There exist oriented graphs~$G$ on $n$ vertices with $\delta^0(G) = \lceil (3n-4)/8 \rceil-1$ without a $1$-factor, see~\cite[Proposition~2]{MR2472138}, and so no $\overset{\rightarrow}{C}_\ell$-factor.
A natural question is whether Theorem~\ref{thm:arboriHC-tiling} holds for some (large) $\ell_0$ independent of~$\eps$. 

\section{Proof of Theorem~\ref{thm:arboriHC-tiling}}

Let $G$ be a digraph.
For a vertex subset~$W$ of $V(G)$, we denote by $G[W]$ the induced subgraph of~$G$ on $W$. 
We write $G \setminus W$ for $G[ V(G) \setminus W]$. 
For $n \in \mathbb{N}$, let $[n] = \{1, \dots, n\}$. 

The key ingredient of Theorem~\ref{thm:arboriHC-tiling} is the following lemma, which states that there exists an equipartition of~$V(G)$ into $V_1, \dots, V_m$ each of size $\ell$, such that $\delta^0(G[V_i])/|V_i|$ is roughly the same as~$\delta^0(G)/|G|$.
The proof of the lemma is based on the proof of~\cite[Lemma~7.4]{BKLO}.
The main difference between~\cite[Lemma~7.4]{BKLO} and the following lemma is that here we are able to control the size of each vertex subset precisely.

\begin{lemma} \label{lma:partition-structure}
There exists $\ell_1 \in \mathbb{N}$ such that for all $\ell \ge \ell_1$ the following holds. 
Let $\delta>0$ and $m, k \in \mathbb{N}$ with $m \le 2^k$. 
Let $G$ be a digraph on $\ell m$ vertices with $\delta^0(G) \ge \delta \ell m $.
Then there exists an equipartition of $V(G)$ into $V_1, \dots, V_m$ such that, for all $i \in [m]$, $|V_i| = \ell$ and 
\begin{align*}
	\delta^0(G[V_i]) \ge \left(\delta -  \ell^{-1/3} \sum_{0 \le j \le k-1} 2^{-j/3} \right) \ell
	\ge \left(\delta - 10\ell^{-1/3}  \right) \ell
	.
\end{align*}
\end{lemma}

Theorem~\ref{thm:arboriHC-tiling} follows immediately from Lemma~\ref{lma:partition-structure} and Theorem~\ref{thm:arboriHC}.

\begin{proof}[Proof of Theorem~\ref{thm:arboriHC-tiling}]
Let $\ell_0(\eps) = \max\{ 20^3\eps^{-3} , n_0, \ell_1 \}$, where $n_0$ and $\ell_1$ are given by Theorem~\ref{thm:arboriHC} and Lemma~\ref{lma:partition-structure}, respectively. 
Let $G,\ell,C_\ell$ be as stated in Theorem~\ref{thm:arboriHC-tiling}.
By Lemma~\ref{lma:partition-structure}, there exists an equipartition of $V(G)$ into $V_1, \dots, V_m$ such that, for all $i \in [m]$,  $|V_i| = \ell$ and 
\begin{align*}
	\delta^0(G[V_i]) \ge \left(\frac38+\eps - 10\ell^{-1/3}  \right) \ell \ge \left(\frac38+\frac{\eps}2 \right)\ell
	.
\end{align*}
By Theorem~\ref{thm:arboriHC}, each $G[V_i]$ contains a spanning copy of $C_\ell$. 
Thus $G$ contains a $C_{\ell}$-factor. 
\end{proof}

Our proof actually proves the slightly stronger result that each cycle on $\ell$ vertices can have a different orientation. 

\subsection{Proof of Lemma~\ref{lma:partition-structure}}
Let $m,n, N \in \mathbb{N}$ with $n,m < N$.
Recall that the hypergeometric distribution with parameters $N, n$ and $m$ is the distribution of the random variable $X$ defined as follows.
Let $S$ be a random subset of $[N]$ of size $n$ and let $X := |S \cap [m]|$.
We use the following simple form of Azuma--Hoeffding's inequality.

\begin{lemma}[see {\cite[Theorem 2.10]{JLR}}] \label{lma:chernoff}
Let $X$ be a hypergeometric distribution with parameters $N,n,m$.
Then $
\mathbb{P}(|X - \mathbb{E}(X)| \geq t) \leq 2e^{-2t^2/n}.
$
\end{lemma}

Let $G$ be a digraph. 
The following proposition is a simple consequence of Lemma~\ref{lma:chernoff}.

\begin{proposition} \label{prop:random-partition}
Let $0 < \alpha \le 1/2$. 
There exists $n_0(\alpha)$ such that the following holds for all $n \ge n_0(\alpha)$ and $m \in \mathbb{N}$ with $\alpha n \le m \le n/2$. 
Let $G$ be a digraph on $n$ vertices with $\delta^0(G) \ge \delta n$ with $0 < \delta <1$.  
Then there exists a vertex subset~$W$ of~$V(G)$ such that $|W| = m$ and 
\begin{align*}
\frac{\delta^0(G[W])}{m} , \frac{\delta^0(G \setminus W)}{n - m}  \ge \delta  -  n^{-1/3}.
\end{align*}
\end{proposition}

\begin{proof}
Consider a random subset $W$ of size~$m$.
Consider any $v \in V(G)$ and $\ast \in \{+,-\}$. 
Let $d^\ast(v,W) = |N^{\ast}(v) \cap W|$.
Note that $d^\ast(v,W)$ has a hypergeometric distribution with parameters $n, m,d^\ast(v)$.
By Lemma~\ref{lma:chernoff} we have that
\begin{align*}
	\mathbb{P} \left( d^*(v,W) \le (\delta  - n^{-1/3})m \right)
	\le 
	\mathbb{P} \left( d^*(v,W) \le \delta m -  n^{-1/3}m \right)  \le  2 e^{ - 2n^{-2/3}m } \le  2 e^{ - 2\alpha n^{1/3} }.
\end{align*}
Similarly, we have 
\begin{align*}
	\mathbb{P} \left( d^*(v,V(G) \setminus W) \le ( \delta - n^{-1/3} )(n- m) \right)  
	\le  2 e^{ - 2n^{-2/3}(n-m) }
	\le  2 e^{ - 2\alpha n^{1/3} }.
\end{align*}
Therefore, together with the union bound, we have
\begin{align*}
\mathbb{P} \left( \frac{\delta^0(G[W])}{m}  , \frac{\delta^0(G \setminus W)}{n - m}  \ge \delta  -  n^{-1/3} \right)
\ge 1  - 2n \cdot 2 e^{ - 2\alpha n^{1/3} } >1/2
\end{align*}
as $n$ is sufficiently large. 
We pick one such~$W$.
\end{proof}

\begin{proof}[Proof of Lemma~\ref{lma:partition-structure}]
Let $\ell_1 = n_0(1/4)$, where $n_0$ is the function given by Proposition~\ref{prop:random-partition}.
Let $G, \ell, \delta,m,k$ be as stated in Lemma~\ref{lma:partition-structure}.
We proceed by induction on~$k$. 
When $k = 1$, we have $m \le 2$. 
If $m=1$, then the lemma holds trivially. 
If $m =2$, then Proposition~\ref{prop:random-partition} with $\alpha = 1/2$ implies the lemma. 
Thus we may assume that $k \ge 2$ and $2^{k-1} < m \le 2^k$.
We may further assume that $\delta >  \ell^{-1/3} \sum_{0 \le j \le k-1} 2^{-j/3}$ or else the lemma holds trivially. 

Let $n = m \ell$ and $n_1 = \lfloor m/2 \rfloor \ell$, so $1/4 \le \lfloor m/2 \rfloor/m = n_1/n \le 1/2$.
By Proposition~\ref{prop:random-partition}, there exists a vertex subset~$W$ of~$V(G)$ such that $|W| = n_1$ and 
\begin{align*}
\frac{\delta^0(G[W])}{|W|} , \frac{\delta^0(G \setminus W)}{|G \setminus W|} \ge \delta  -  n^{-1/3} \ge \delta  - \left( 2^{k-1} \ell \right)^{-1/3},
\end{align*} 
where the last inequality holds as $n = m \ell > 2^{k-1} \ell$.
Note that $|G[W]|, |G \setminus W| \le 2^{k-1} \ell$.
Apply our induction hypothesis to both $G[W]$ and $G \setminus W$ and obtain an equipartition of $V(G)$ into $V_1, \dots, V_m$ such that, for all $i \in [m]$, $|V_i| = \ell$ and 
\begin{align*}
\delta^0(G[V_i]) & \ge \left(  (\delta -  2^{-(k-1)/3} \ell^{-1/3}) -  \ell^{-1/3} \sum_{0 \le j \le k-2} 2^{-j/3} \right) \ell
= \left( \delta -  \ell^{-1/3} \sum_{0 \le j \le k-1} 2^{-j/3} \right) \ell \\
	& \ge \left( \delta -  \ell^{-1/3} \sum_{j \ge 0} 2^{-j/3} \right) \ell 
	= \left( \delta -  \ell^{-1/3} \frac{1}{1-2^{-1/3}} \right) \ell
	\ge \left( \delta - 10 \ell^{-1/3} \right) \ell
\end{align*}
as required. 
\end{proof}


\section{Remarks}

Here, we discuss other results that can be obtained using a similar argument. 

\subsection{Powers of directed cycles}

The $k$th power of a digraph~$G$ is the digraph obtained from~$G$ by adding a directed edge from~$x$ to~$y$ if there is a directed path of length at most $k$ from $x$
to $y$ on~$G$. 
Let ${C}^{k}_\ell$ be the $k$th power of a directed cycle on $\ell$ vertices. 
Debiasio, Han, Molla, Piga, Treglown and the author~\cite[Theorem~1.5]{debiasio2024powers} proved a minimum semi-degree threshold for the $k$th power of a directed Hamilton cycle in oriented graphs. 
Our argument, namely Lemma~\ref{lma:partition-structure}, gives the following result.

\begin{corollary} 
For all $k \ge 2$, there exists $\ell_0 = \ell_0(k)$ such that for all $\ell \ge \ell_0$, any oriented graph~$G$ on $ \ell m$ vertices with minimum semi-degree $\delta^0(G)\geq (1/2- {10^{-6001k}} ) \ell m$ contains ${C}^{k}_\ell$-factor.
\end{corollary}

\subsection{Other minimum degree conditions}
Let $G$ be a digraph. 
Recall that the \emph{minimum total degree}~$\delta(G)$ of $G$ is the minimum $\min \{d^+(x) + d^-(x) : x \in V(G)\}$.
One can prove a corresponding version of Lemma~\ref{lma:partition-structure} for minimum total degree.\footnote{Indeed, the proof of Proposition~\ref{prop:random-partition} yields a partition~$W$ such that for all $\ast \in \{+,-\}$ and all $v \in V(G)$, we have $d^*(v,W)/m , d^*(v,V(G) \setminus W)/(n-m) \ge d^*(v)/n - n^{-1/3}$.
The rest of the argument follows similarly to the proof of Lemma~\ref{lma:partition-structure}.}
Thus, we obtain the following results, which follow from \cite[Theorem 1.4]{debiasio2024powers} and a result of DeBiasio and Treglown~\cite{debiasio2025arbitrary}, respectively.

\begin{corollary} 
For all $\eps >0$, there exists $\ell_0 = \ell_0(\eps)$ such that for all $\ell \ge \ell_0$, any digraph~$G$ on $ \ell m$ vertices with minimum total degree $ \delta(G)\geq (8/5+\eps ) \ell m$ contains a ${C}^{2}_\ell$-factor.
\end{corollary}

\begin{corollary} 
For all $\eps >0$, there exists $\ell_0 = \ell_0(\eps)$ such that for all $\ell \ge \ell_0$ and any oriented cycle~$C_\ell$ (except the directed cycle) on $\ell$ vertices, any digraph~$G$ on $ \ell m$ vertices with minimum total degree $\delta(G) \ge  (1 + \eps) \ell m$ contains a $C_\ell$-factor.
\end{corollary}

One can also prove a $k$-uniform hypergraph version of Lemma~\ref{lma:partition-structure} for the minimum $d$-degree. 
However, we are unaware of any immediate improvement to existing $H$-factor results using this method. 

\subsection{Limitation of the method}

The minimum degree threshold for finding an $H$-factor obtained using the method presented is unlikely to be asymptotically sharp.
This is because $\eps$ is a function of~$|H|$. 
Our proof gives $\eps = O(|H|^{-1/3})$.
In some cases, $\eps$ cannot be removed and depends on the `parity' of $|H|$.
For example, suppose that we apply this method to obtain bound for finding a $C_\ell$-factor in graphs. 
Using Dirac's theorem, this would yield a bound of $\delta(G) \ge (1/2 + \eps_\ell) n $, where $\eps_\ell$ tends to zero as $\ell$ tends to infinity.
However, when $\ell$ is odd, $\eps_\ell$ cannot be removed.
In particular, $\delta(G)  \ge (1/2 + 1/(2\ell))n +c$ is the correct threshold for $C_\ell$-factor up to an additive constant~$c$, see~\cite{MR2506388}.

\bibliographystyle{abbrv}
\bibliography{reference}

\begin{thebibliography}{10}

\bibitem{BKLO}
B.~Barber, D.~K\"uhn, A.~Lo, and D.~Osthus.
\newblock Edge-decompositions of graphs with high minimum degree.
\newblock {\em Adv. Math.}, 288:337--385, 2016.

\bibitem{debiasio2024powers}
L.~DeBiasio, J.~Han, A.~Lo, T.~Molla, S.~Piga, and A.~Treglown.
\newblock Powers of {H}amilton cycles in oriented and directed graphs.
\newblock {\em Combin. Probab. Comput.}, 35(1):101--133, 2026.

\bibitem{debiasio2025arbitrary}
L.~DeBiasio and A.~Treglown.
\newblock Arbitrary orientations of {H}amilton cycles in directed graphs of
  large minimum degree.
\newblock {\em arXiv preprint arXiv:2505.09793}, 2025.

\bibitem{Dirac}
G.~A. Dirac.
\newblock Some theorems on abstract graphs.
\newblock {\em Proc. London Math. Soc. (3)}, 2:69--81, 1952.

\bibitem{HajnalSzemeredi}
A.~Hajnal and E.~Szemer\'edi.
\newblock Proof of a conjecture of {P}. {E}rd{\H{o}}s.
\newblock In {\em Combinatorial theory and its applications, {I}-{III} ({P}roc.
  {C}olloq., {B}alatonf\"ured, 1969)}, volume~4 of {\em Colloq. Math. Soc.
  J\'anos Bolyai}, pages 601--623. North-Holland, Amsterdam-London, 1970.

\bibitem{JLR}
S.~Janson, T.~{\L}uczak, and A.~Ruci\'nski.
\newblock {\em Random graphs}.
\newblock Wiley-Interscience Series in Discrete Mathematics and Optimization.
  Wiley-Interscience, New York, 2000.

\bibitem{MR2472138}
P.~Keevash, D.~K\"uhn, and D.~Osthus.
\newblock An exact minimum degree condition for {H}amilton cycles in oriented
  graphs.
\newblock {\em J. Lond. Math. Soc. (2)}, 79(1):144--166, 2009.

\bibitem{KeevashSudakov}
P.~Keevash and B.~Sudakov.
\newblock Triangle packings and 1-factors in oriented graphs.
\newblock {\em J. Combin. Theory Ser. B}, 99(4):709--727, 2009.

\bibitem{Kelly_arborientation}
L.~Kelly.
\newblock Arbitrary orientations of {H}amilton cycles in oriented graphs.
\newblock {\em Electron. J. Combin.}, 18(1):Paper 186, 25, 2011.

\bibitem{MR2454564}
L.~Kelly, D.~K\"uhn, and D.~Osthus.
\newblock A {D}irac-type result on {H}amilton cycles in oriented graphs.
\newblock {\em Combin. Probab. Comput.}, 17(5):689--709, 2008.

\bibitem{MR2588541}
D.~K\"uhn and D.~Osthus.
\newblock Embedding large subgraphs into dense graphs.
\newblock In {\em Surveys in combinatorics 2009}, volume 365 of {\em London
  Math. Soc. Lecture Note Ser.}, pages 137--167. Cambridge Univ. Press,
  Cambridge, 2009.

\bibitem{MR2506388}
D.~K\"uhn and D.~Osthus.
\newblock The minimum degree threshold for perfect graph packings.
\newblock {\em Combinatorica}, 29(1):65--107, 2009.

\bibitem{MR2889513}
D.~K\"uhn and D.~Osthus.
\newblock A survey on {H}amilton cycles in directed graphs.
\newblock {\em European J. Combin.}, 33(5):750--766, 2012.

\bibitem{MR4025428}
L.~Li and T.~Molla.
\newblock Cyclic triangle factors in regular tournaments.
\newblock {\em Electron. J. Combin.}, 26(4):Paper No. 4.24, 23, 2019.

\bibitem{wang2025arbitrary}
G.~Wang, Y.~Wang, and Z.~Zhang.
\newblock Arbitrary orientations of cycles in oriented graphs.
\newblock {\em arXiv preprint arXiv:2504.09794}, 2025.

\bibitem{MR4763245}
Z.~Wang, J.~Yan, and J.~Zhang.
\newblock Cycle-factors in oriented graphs.
\newblock {\em J. Graph Theory}, 106(4):947--975, 2024.

\bibitem{MR3526407}
Y.~Zhao.
\newblock Recent advances on {D}irac-type problems for hypergraphs.
\newblock In {\em Recent trends in combinatorics}, volume 159 of {\em IMA Vol.
  Math. Appl.}, pages 145--165. Springer, [Cham], 2016.

\end{thebibliography}



\end{document}